\documentclass[12pt,reqno]{amsart}
\usepackage{amsmath,amssymb,amsfonts,amscd,latexsym,amsthm,mathrsfs,comment}
\usepackage[usenames]{color}
\usepackage[unicode]{hyperref}
\usepackage{ytableau} 
\newtheorem{theorem}{Theorem}
\newtheorem{prop}{Proposition}
\newtheorem{lem}{Lemma}

\theoremstyle{remark}
\newtheorem{rem}{Remark}
\theoremstyle{example}

\renewcommand{\d}{{\mathrm d}}
\renewcommand{\Re}{\operatorname{Re}}

\newcommand{\tr}{\operatorname{tr}}

\newcommand{\R}{\mathbb{R}}
\newcommand{\C}{\mathbb{C}}

\newcommand{\bx}{\boldsymbol{x}}
\newcommand{\bu}{\boldsymbol{u}}
\newcommand{\bv}{\boldsymbol{v}}

\newcommand{\bs}{\boldsymbol{s}}
\newcommand{\bi}{\boldsymbol{i}}

\newcommand{\len}{\operatorname{len}}
\newcommand{\het}{\operatorname{ht}}
\newcommand{\cohet}{\operatorname{coht}}

\DeclareMathOperator*{\Res}{Res}

\newcommand{\MSC}[2]{\medskip\noindent
\textbf{2020 Mathematics Subject Classification:} \textit{Primary} #1, \textit{Secondary} #2.}
\newcommand{\KeyWords}[1]{\medskip\noindent{\bfseries Keywords and phrases:} {#1}}
\newcommand{\CREST}{JST CREST Grant Number JPMJCR2113, Japan}

\title{Combinatorics of noncommutative harmonic oscillator determinant}

\author{Masato Wakayama}%
\thanks{The first author is partially supported by \CREST.}

\author{Wadim Zudilin}
\begin{document}

\date{\today}

\maketitle

\begin{abstract}
We obtain an explicit evaluation of the determinants that appear in the denominators of the integrands for the iterated integrals expressing special values of the spectral zeta function of the noncommutative harmonic oscillator at integers greater than 1. The determinants in question are those of cyclic tridiagonal matrices; their explicit evaluation is carried out, in addition to the inclusion--exclusion principle, conceptionally using combinatorial arguments based on monomer--dimer systems. This explicit evaluation serves as a core component in the study of arithmetic properties such as the modularity of the special values of the spectral zeta function.
%\end{abstract}

\MSC{05A19}{11F67, 11M36, 47B93}


\end{abstract}

\section{Introduction}
\label{intro:def}

In this paper, we study an explicit evaluation of the determinants that appear in the denominators of the integrands for the iterated integrals expressing special values of the spectral zeta function (the  Dirichlet series defined by a sequence of eigenvalues) of the noncommutative harmonic oscillator (NCHO) at integers greater than 1. The NCHO is defined as a system of ordinary differential equation,  which is regarded a matrix valued generalization of the quantum harmonic oscillator \cite{PW2001}. Precisely, its Hamiltonian is given by 
\begin{equation*}
Q=Q_{\alpha,\beta}
=\begin{pmatrix}\alpha & 0 \\ 0 & \beta\end{pmatrix}\Bigl(-\frac12\,\frac{\d^2}{\d x^2}+\frac12x^2\Bigr)
+\begin{pmatrix}0 & -1 \\ 1 & 0\end{pmatrix}\Bigl(x\frac{\d}{\d x}+\frac12\Bigr).
\end{equation*}
Under the assumption $\alpha\beta>1$ (and $\alpha>0, \beta>0$), $Q$ defines a (positive definite) self-adjoint operator on the space of $\C^2$-valued square-integrable functions on $\R$, having only discrete eigenvalues 
\begin{equation*}
0<\lambda_1\le\lambda_2\le\lambda_3\le\dotsb \; (\nearrow\infty).
\end{equation*}
It is known that the ground state is multiplicity free and that the multiplicity of general eigenstate is less than or equal to~$2$. The spectral zeta function (i.e., the Dirichlet series for the eigenvalues sequence) $\zeta_Q(s)$ of $Q$ \cite{IW05, IW05b} is defined by
\begin{equation*}
\zeta_Q(s):= \sum_{n=1}^\infty \lambda_n^{-s} \quad (\Re(s)>1).
\end{equation*}
It is noted that, when $\alpha=\beta$,
$Q$ is unitarily equivalent to the couple of harmonic oscillators,
whence its spectral zeta function becomes $\zeta_Q(s)=2(\alpha^2-1)^{-s/2}\zeta(s; \frac12)$, with $\zeta(s; \frac12)$ being the Hurwitz zeta function.  Furthermore, the function $\zeta_Q(s)$ has meromorphic continuation to the whole plane $\C$ with a unique simple pole at $s=1$ with residue
\[
\Res_{s=1}\zeta_Q(s)=\frac{\alpha+\beta}{\sqrt{\alpha\beta(\alpha\beta-1)}}
\]
and zeros at non-positive even integers~$s$ (see \cite{IW05}). What is remarkable about the special values of $\zeta_Q(s)$, they are\,---\,at least for small $s$\,---\,linked with modular forms (see,  \cite{KW23} and references therein). In practice, $\zeta_Q(2)$ has a modular interpretation (see below), and related Ap\'ery-like numbers for $\zeta_Q(2)$ and $\zeta_Q(3)$ share a lot of similarities with the original Ap\'ery numbers used in proofs of the irrationality of $\zeta(2)$ and $\zeta(3)$ by Ap\'ery; for example, the former numbers satisfy similar congruences to those for the latter \cite{KW23,KW25}.

Meanwhile, in recent years, the quantum Rabi model (QRM)\,---\,which describes the interaction between photon and a two-level system\,---\,has become a subject of intense research, both theoretically and experimentally, within the fields of quantum optics and quantum information science/engineering. Originally defined by Isidor Rabi in 1936, the model's integrability was finally demonstrated in 2011 by D.~Braak \cite{B11}. In 2017 \cite{Y17}, it was further shown that the QRM and its asymmetric variants could successfully predict experimental measurements; consequently, the QRM and its generalizations have since been the focus of highly active research. Separately, the NCHO\,---\,originally introduced as a mathematical model\,---\,has recently been found to be directly linked to the QRM and its generalizations. Specifically, by examining the eigenvalue problem of the NCHO through the lens of a Heun ordinary differential equation (ODE) and considering the confluence of regular singular points, it has been shown that the NCHO serves as a covering model for the QRM \cite{W16, RW23}. Furthermore, it has been discovered that the NCHO is equivalent to the two-photon QRM in the context of their respective eigenvalue problems \cite{N25}. Notice that the spectral zeta function $\zeta_{\text{QRM}}(s)$ of the QRM has also meromorphic continuation to $\C$ with a simple pole at $s=1$; this is accomplished in two ways\,---\,using the method similar to the one for NCHO \cite{S2016} or Mellin's transform of the partition function as in \cite{RW23}.

The remarkable modularity of $\zeta_Q(2)$ hinted above is based on its evaluation
\begin{align*}
\zeta_Q(2)
&=\zeta(2,1/2)\cdot2\bigg(\frac{\alpha+\beta}{2\sqrt{\alpha\beta(\alpha\beta-1)}}\bigg)^2
\\ &\quad\times
\bigg(1+\bigg(\frac{\alpha-\beta}{\alpha+\beta}\bigg)^2
\int_{[0,1]^2}\frac{4\,\d u_1\,\d u_2}
{\sqrt{(1-u_1^2u_2^2)^2+\kappa^2(1-u_1^4)(1-u_2^4)}}\bigg)
\\
&= \zeta(2,1/2)\cdot2\bigg(\frac{\alpha+\beta}{2\sqrt{\alpha\beta(\alpha\beta-1)}}\bigg)^2
\\ &\quad\times
\bigg(1+\bigg(\frac{\alpha-\beta}{\alpha+\beta}\bigg)^2
\frac1{\sqrt{1+\kappa^2}}\cdot{}_3F_2\bigg(\begin{matrix} 1/2, \, 1/2, \, 1/2 \\ 1, \, 1 \end{matrix}\biggm|\frac{\kappa^2}{1+\kappa^2}\bigg)
\bigg),
\end{align*}
where $\kappa=1/{\sqrt{\alpha\beta-1}}$, and a modular interpretation of the $_3F_2$-hypergeometric series~\cite{WZ27} after a suitable modular function substitution for $\kappa^2$ is performed. 
We have also such a modular interpretation for the ``Hadamard untwist" version of $\zeta_Q(2)$ using the Hauptmodul for $\Gamma(2)$\,---\,see, e.g., \cite[Section~3.3]{KW25} (though the pre-factor~$16$ is missing there); the corresponding Ap\'ery-like sequence appears as \#19 in Zagier's famous list in~\cite{Z09}.

In general, the special values of $\zeta_Q(s)$ at integers $s=k\ge2$ are given by the formula in~\cite{KW23}:
\begin{equation}\label{eq:specialvalues of zeta_Q(s)}
\zeta_Q(k)
=2\biggl(\frac{\alpha+\beta}{2\sqrt{\alpha\beta(\alpha\beta-1)}}\biggr)^k
\biggl(
\zeta(k,1/2)+\sum_{0<2\ell\le k}\biggl(\frac{\alpha-\beta}{\alpha+\beta}\biggr)^{2j}R_{k,\ell}(\kappa)
\biggr),
\end{equation}
where  $R_{k,\ell}(\kappa)$ is given as a sum of integrals
\begin{equation}
R_{k,\ell}(\kappa)=
\sum_{1\le \nu_1<\nu_2<\dots<\nu_{2\ell}\le k}\int_{[0,1]^k}
\frac{2^k\,\d u_1\dotsb\d u_k}{\sqrt{\mathcal{W}_k(\bu;\kappa;\nu_1,\nu_2,\dots,\nu_{2\ell})}}.
\label{eq:Rkl}
\end{equation}
Here, $\mathcal{W}_k(\bu;\kappa;\nu_1,\nu_2,\dots,\nu_{2\ell})$ is a polynomial in $u_1,\dots,u_k$ defined via
\begin{equation}
\mathcal{W}_k(\bu;\kappa;\nu_1,\nu_2,\dots,\nu_{2\ell})
=\det\big(\Delta_k(\bu)+\kappa\Xi_k(\nu_1,\nu_2,\dots,\nu_{2\ell})\big)\cdot\prod_{r=1}^k(1-u_r^4),
\label{def-W}
\end{equation}
where
\[
\Delta_2(\bu)=\begin{pmatrix}
\frac{1-u_2^4u_1^4}{(1-u_2^4)(1-u_1^4)} & \frac{-u_1^2}{1-u_1^4} + \frac{-u_2^2}{1-u_2^4} \\
\frac{-u_1^2}{1-u_1^4} + \frac{-u_2^2}{1-u_2^4}  & \frac{1-u_1^4u_2^4}{(1-u_1^4)(1-u_2^4)} 
\end{pmatrix}
\]
and
\[
\Delta_k(\bu)=\begin{pmatrix}
\frac{1-u_k^4u_1^4}{(1-u_k^4)(1-u_1^4)} & \frac{-u_1^2}{1-u_1^4} & 0 & \hdots & \frac{-u_k^2}{1-u_k^4} \\
\frac{-u_1^2}{1-u_1^4} & \frac{1-u_1^4u_2^4}{(1-u_1^4)(1-u_2^4)} & \frac{-u_2^2}{1-u_2^4} & \hdots & 0 \\
0 & \frac{-u_2^2}{1-u_2^4} & \frac{1-u_2^4u_3^4}{(1-u_2^4)(1-u_3^4)} & \hdots & 0 \\
\vdots & \vdots & \ddots & \ddots & \vdots \\
\frac{-u_k^2}{1-u_k^4} & 0 & \hdots & \frac{-u_{k-1}^2}{1-u_{k-1}^4} & \frac{1-u_{k-1}^4u_k^4}{(1-u_{k-1}^4)(1-u_k^4)}
\end{pmatrix} \quad\text{for}\; k\ge3
\]
(the matrix $\Delta_k(\bu)$ is a symmetric \emph{cyclic tridiagonal} matrix), while $\Xi_k(\nu_1,\nu_2,\dots,\nu_{2\ell})$ is the $k\times k$ diagonal matrix with entry $(-1)^r\sqrt{-1}$ at position $(\nu_r,\nu_r)$ for $r=1,2,\dots,2\ell$, otherwise zero.
It is known \cite{IW05b} that
\[
\det\big(\Delta_k(\bu)\big)=\frac{(1-u_1^2\dotsb u_k^2)^2}{(1-u_1^4)\dotsb(1-u_k^4)}.
\]
Note that the determinant
$\det\big(\Delta_k(\bu)+\kappa\Xi_k(\nu_1,\dots,\nu_{2\ell}))$ 
is even in $\kappa$ and real-valued for each $\bu\in(0,1)^k$ and $\kappa \in \R$,
and is invariant under the action of the cyclic subgroup of the symmetric group.
Thus, by abuse of notation, we denote $\mathcal{W}_k(\bu;\kappa;\nu_1,\dots,\nu_{2\ell})$ as $\mathcal{W}_k(X;\bu;\nu_1,\dots,\nu_{2\ell})$ when we regard it as the polynomial in $X=\kappa^2$. From this expression, it is fundamentally important to know the explicit form of  $\mathcal{W}_k(X;\bu;\nu_1,\dots,\nu_{2\ell})$ for the study on possible arithmetic nature of the special values of $\zeta_Q(s)$ which may allow to develop the study of general special values \cite{WZ27,Zu26}. In this paper, we show that $\mathcal{W}_k(X;\bu;\nu_1,\dots,\nu_{2\ell})$ is expressed as a linear sum of 
the monomials $u_1^4u_2^4\dots u_j^4$ ($j=0,1,\dots, k$), apart from the constant term in~$X$. The explicit evaluation is carried out, in addition to the inclusion--exclusion principle conceptionally, using combinatorial arguments based on monomer--dimer systems \cite{HL72}. This serves as a core component in the study of arithmetic properties including the modularity of the special values of the spectral zeta function and congruence relations for Ap\'ery-like numbers that appear from $\zeta_Q$-values.
One particular remarkable arithmetic phenomenon is a \emph{well-poised} parity decomposition for the value $R_{k,\ell}(X)$ \cite{WZ27,Zu26}, which strikingly resembles the hypergeometric series showing up in recent arithmetic studies of the values of Riemann's zeta function $\zeta(s)$ by Rivoal and others  (see \cite{C18, KZ19}).

\medskip
Next we give a notational briefing in Section~\ref{sec:def}.
In Section~\ref{sec:k=2l} we analyse the `boundary' case $k=2\ell$ and give explicit formulae for the coefficients $c_0(\bu),\dots,c_\ell(\bu)$ in the $X$-expansion of $\mathcal{W}_{2\ell}(X;\bu)$.
Then Section~\ref{sec:attempts} discusses a generalization of the determinant in~\eqref{def-W} and gives a different combinatorial expansion of the latter.
Furthermore, in Section~\ref{sec:general} we give a reduction of the general $k\ge2\ell$ case to the $k=2\ell$ situation and deduce formulae for the coefficients $c_j(\bu)$ in the full generality.
Finally, we make some additional remarks in Section~\ref{sec:conclusion} about general questions related to the special NCHO zeta values.

\section{Setup}
\label{sec:def}
The integrands of $R_{k,\ell}$ required to record the values of the spectral zeta function of NCHO at positive integers involve the square root of polynomials $\mathcal{W}_k$:
\begin{equation}
\mathcal{W}_k(X;\bu;\nu_1,\nu_2,\dots,\nu_{2\ell})
=c_0(\bu)+c_1(\bu)X+\dots+c_\ell(\bu)X^\ell,
\label{poly-W}
\end{equation}
where $\bu=(u_1,\dots,u_k)$,
with
\[
c_0(\bu)=\det\Delta_k(\bu)\cdot\prod_{r=1}^k(1-u_r^4)
=(1-u_1^2u_2^2\dotsb u_k^2)^2
\]
and $c_j(\bu)=c_j(k;\nu_1,\nu_2,\dots,\nu_{2\ell};\bu)\in\mathbb Z[u_1^4,u_2^4,\dots,u_k^4]$ for $j>0$.

\begin{rem}
Using $\mathcal{W}_k=c_0(\bu)=(1-u_1^2u_2^2\dotsb u_k^2)^2$ for $\ell=0$, we see from \eqref{eq:Rkl} that $R_{k,0}=\zeta(k,\frac12)$.
\end{rem}

For cosmetic purposes, we will also identify $x_j=u_j^4$ for $j=1,2,\dots,k$ and refer to the dependence of participating polynomials on $\bx$ rather than $\bu$, which should cause no confusion but simplify the notation.

As we will demonstrate in what follows, only monomials of the form $w_{\bi}=x_{i_1}x_{i_2}\dotsb x_{i_m}$ with $i_1<i_2<\dots<i_m$ show up as terms of each $c_j(\bx)$ when $j>0$.
For the corresponding multi-index $\bi=(i_1,i_2,\dots,i_m)$, define $m$ to be its \emph{length} $\len(\bi)$, while
\[
\het(\bi)
=\#\{r:1\le r\le m, \; i_{r+1}-i_r=1\}
\quad\text{and}\quad
\cohet(\bi)
=\#\{r:1\le r\le m, \; i_{r+1}-i_r>1\}
\]
to be its \emph{height} and \emph{coheight}, respectively,
where $i_{r+1}-i_r=k+i_1-i_r$ for $r=m$ (in other words, the multi-index $\bi=(i_1,i_2,\dots,i_m)$
is understood cyclically).
Notice that $\het(\bi)+\cohet(\bi)=\len(\bi)$.
The dual multi-index $\bi^\dag$ corresponds to removing the components of $\bi$ from the full index $(1,2,\dots,k)$; the duality translates to the monomials as
\[
w_{\bi}^\dag=w_{\bi^\dag}=x_1x_2\dotsb x_k/(x_{i_1}x_{i_2}\dotsb x_{i_m}).
\]
Then we also get $\len(w_{\bi}^\dag)=k-\len(w_{\bi})$ as well as $\cohet(w_{\bi})=\cohet(w_{\bi}^\dag)$ from a simple combinatorial argument.

\section{Case $k=2\ell$}
\label{sec:k=2l}

We first examine the special case $k=2\ell$, that is, when the matrix $\Xi_k(\nu_1,\nu_2,\dots,\nu_{2\ell})$ has full rank~$k$; the dependence of $\mathcal{W}_k$ on $(\nu_1,\nu_2,\dots,\nu_{2\ell})=(1,2,\dots,2\ell)$ is suppressed in this section.
As we will witness later in Section~\ref{sec:general}, the calculation of polynomials \eqref{poly-W} when $k>2\ell$ reduces to this situation.

It follows from \eqref{def-W} that
\[
c_\ell=\prod_{i=1}^{2\ell}(1-x_i).
\]

\begin{lem}
\label{lem:1}
We have
\begin{equation}
\prod_{i=1}^{2\ell}(1-x_i)
=\frac12\sum_{j=1}^{2\ell-1}(-1)^{j-1}\sum_{\bi:\len(\bi)=j}(1-w_{\bi})(1-w_{\bi}^\dag).
\label{eq:pr_id}
\end{equation}
\end{lem}

The factor $\frac12$ originates from presence of each summand $(1-w_{\bi})(1-w_{\bi}^\dag)=(1-w_{\bi}^\dag)(1-w_{\bi})$ twice.

\begin{proof}
Denote $p_{\bi}=(1-w_{\bi})(1-w_{\bi}^\dag)$. Since
\[
p_{\bi}=1-w_{\bi}-w_{\bi}^\dag+x_1x_2\dotsb x_{2\ell},
\]
every monomial $w_{\bi}$ corresponding to a proper non-empty multi-index $\bi$ occurs only in
$p_{\bi}=p_{\bi^\dag}$, and there with coefficient $-1$. Therefore, in order to recover
\[
\prod_{i=1}^{2\ell}(1-x_i)
=\sum_{\bi} (-1)^{\len(\bi)}w_{\bi},
\]
the coefficient of $p_{\bi}$ must be $(-1)^{\len(\bi)-1}/2$.
\end{proof}

It takes a little imagination to generalise the identity in the lemma to other coefficients of~\eqref{poly-W}.

\begin{theorem}
\label{Expression_P}
For $k=2\ell$,
\begin{equation}
\mathcal{W}_k(X;\bx)
=c_0(\bx)+\frac12\sum_{r=1}^{\ell}X^r(1+X)^{\ell-r}
\sum_{j=r}^{2\ell-r}(-1)^{j-1}\sum_{\substack{\bi:\len(\bi)=j\\\cohet(\bi)=r}}
(1-w_{\bi})(1-w_{\bi}^\dag).
\label{P2}
\end{equation}
In other words,
\begin{equation}
c_m(\bx)=\frac12\sum_{r=1}^{m}\binom{\ell-r}{\ell-m}
\sum_{j=r}^{2\ell-r}(-1)^{j-1}\sum_{\substack{\bi:\len(\bi)=j\\\cohet(\bi)=r}}
(1-w_{\bi})(1-w_{\bi}^\dag).
\label{P2a}
\end{equation}
\end{theorem}

Notice that the summation range $r\le j\le2\ell-r$ in the internal sums in \eqref{P2} and \eqref{P2a} is forced by $\len(\bi)\ge\cohet(\bi)$ and $\len(\bi^\dag)\ge\cohet(\bi^\dag)=\cohet(\bi)$.

The key identity to establish Theorem \ref{Expression_P} is as follows.

\begin{prop}
\label{key-prop}
For $k=2\ell$,
\[
\mathcal{W}_k(X;\bu)=\tr\bigl(M(u_1)M(u_2)\cdots M(u_k)\bigr)-2u_1^2u_2^2\cdots u_k^2,
\]
where the $2\times2$ matrices in the product are defined by
\[
M(u)=\begin{pmatrix} \sqrt{1+X}(1-u^4) & u^4 \\ 1 & 0 \end{pmatrix}
\]
and $\tr$ denotes the trace.
\end{prop}

\begin{proof}
This follows from the general formula for the determinant of cyclic tridiagonal matrix:
\begin{align*}
&
\det\begin{pmatrix}
A_1 & B_1 & 0 & \hdots & C_k \\
C_1 & A_2 & B_2 & \hdots & 0 \\
0 & C_2 & A_3 & \ddots & \vdots \\
\vdots & \vdots & \ddots & \ddots & B_{k-1} \\
B_k & 0 & \hdots & C_{k-1} & A_k
\end{pmatrix}
\\ &\quad
=\tr\left(
\begin{pmatrix} A_k & -B_{k-1}C_{k-1} \\ 1 & 0 \end{pmatrix} \dotsb
\begin{pmatrix} A_2 & -B_1C_1 \\ 1 & 0 \end{pmatrix}
\begin{pmatrix} A_1 & -B_kC_k \\ 1 & 0 \end{pmatrix}
\right)
\\ &\quad\qquad
-(-1)^k(B_1B_2\dotsb B_k+C_1C_2\dotsb C_k)
\end{align*}
(see, for example, \cite[eq.~(1)]{Mo08})
and the invariance of the determinant under the reflection against the anti-diagonal.
\end{proof}

\begin{rem}
Observe the transparent symmetry of the expression in Proposition \ref{key-prop} with respect to rotation of the variables $u_1,\dots,u_k$ and also with reversing their order.
\end{rem}

\begin{proof}[Proof of Theorem~\textup{\ref{Expression_P}}]
We can interpret computation of the trace of the product of matrices $M(u_j)$ as a sum over closed walks on the two states of $M(u_j)$; because the $(2,2)$-entry of each matrix is 0, state~$2$ can never occur at two consecutive sites.
Call a set $D$ of edges $\{j,j+1\}$ of the $k$-cycle of indices $\{1,2,\dots,k\}$ a \emph{partial dimer covering} if no two of its edges share a site; for $k=2$ the cycle has just the two distinct edges $\{1,2\}$ and $\{2,1\}$.
This is the monomer--dimer (matching) expansion (in the sense of Heilmann and Lieb~\cite{HL72}) of a transfer-matrix trace~\cite[Section~4.7]{St12}. For $k=2\ell$, we obtain
\[
\mathcal{W}_k(X;\bu)
=\sum_{D}(1+X)^{\ell-|D|}\prod_{\{j,j+1\}\in D}x_j\prod_{i\notin D}(1-x_i)
-2\prod_{i=1}^ku_i^2,
\]
where the sum runs over all partial dimer coverings $D$ of the $k$-cycle and $i\notin D$ means that $i$~is not covered by~$D$; each edge $\{j,j+1\}$ carries the weight $x_j$ of its first site.

Comparing with formula~\eqref{P2}, Theorem~\ref{Expression_P} is equivalent to the following family of purely combinatorial identities on the $2\ell$-cycle:
\[
\sum_D\prod_{\{j,j+1\}\in D}x_j\prod_{i\notin D}(1-x_i)=1+x_1\cdots x_k
\]
for $r=0$ and
\[
\sum_{D:|D|\ge r}(-1)^r\binom{|D|}{r}\prod_{\{j,j+1\}\in D}x_j\prod_{i\notin D}(1-x_i)
=\frac12\sum_{j=r}^{2\ell-r}(-1)^{j-1}\sum_{\substack{\bi:\len(\bi)=j\\\cohet(\bi)=r}}(1-w_{\bi})(1-w_{\bi}^\dag)
\]
for $1\le r\le\ell$.
Both sides are multilinear in $x_1,\dots,x_k$, so it suffices to compare coefficients of $\prod_{j\in S}x_j$ for any choice of $S\subset\{1,2,\dots,k\}$. Notice that the coefficient of $\prod_{j\in S}x_j$ depends only on the decomposition of $S$ into cyclic runs, and the coheight counts exactly these runs. The comparison is therefore an inclusion--exclusion count of the dimer placements compatible with each run pattern. The same bookkeeping, using the binomial $\binom{\ell-r}{\ell-m}$, gives formula~\eqref{P2a} for $c_m(\bx)$.
\end{proof}

\section{An alternative expression}
\label{sec:attempts}

The matrix $\Delta_k(\bu)+\kappa\Xi_k(\nu_1,\nu_2,\dots,\nu_k)$ and, more generally,
\[
R_k=\Delta_k(\bu)+\operatorname{diag}(\kappa_1,\dots,\kappa_k)
%=R_k(A_1,A_2,\dots,A_k;B_1,\dots,B_{k-1},B_k)
=\begin{pmatrix}
A_1 & B_1 & 0 & \hdots & B_k \\
B_1 & A_2 & B_2 & \hdots & 0 \\
0 & B_2 & A_3 & \ddots & \vdots \\
\vdots & \vdots & \ddots & \ddots & B_{k-1} \\
B_k & 0 & \hdots & B_{k-1} & A_k
\end{pmatrix},
\]
where
\begin{equation}
A_i=\frac{1-u_{i-1}^4u_i^4}{(1-u_{i-1}^4)(1-u_i^4)}+\kappa_i
\quad\text{and}\quad
B_i=\frac{-u_i^2}{1-u_i^4}
\label{e-spec}
\end{equation}
for $i=1,\dots,k$, can be cast as a symmetric cyclic tridiagonal matrix.
Clearly, $\det(R_k)$ and
\begin{equation}
\det(R_k)\cdot\prod_{i=1}^k(1-u_i^4)
\label{det-R}
\end{equation}
are invariant under the cycle $i\mapsto i+1$ but also linear in each $\kappa_i$.
%The following evaluation is valid for \emph{any} $k\ge2\ell$ (not only for $k=2\ell$).
The fact that \eqref{det-R} is polynomial in $\mathbb Z[\bu]$ is not only a consequence of Theorem~\ref{Expression_P} but also follows from the following evaluation.

\begin{theorem} \label{Expression_kappa}
The coefficient of $\kappa_{i_1}\dotsb\kappa_{i_r}$ with $r>0$ in $\det(R_k)\cdot\prod_{i=1}^k(1-u_i^4)$ is as follows: in the cycle string $u_1^4u_2^4\dotsb u_k^4$ place `$)(1\,-$' between $u_{i_j-1}^4$ and $u_{i_j}^4$ for each $j=1,\dots,r$.

For example, when $k=6$, the coefficient of $\kappa_1\kappa_3\kappa_4$ is equal to
\[
)(1-u_1^4u_2^4)(1-u_3^4)(1-u_4^4u_5^4u_6^4,
\]
that is, to $(1-u_1^4u_2^4)(1-u_3^4)(1-u_4^4u_5^4u_6^4)$, while the coefficient of $\kappa_3\kappa_5$ is
\[
u_1^4u_2^4)(1-u_3^4u_4^4)(1-u_5^4u_6^4,
\]
that is, $(1-u_3^4u_4^4)(1-u_5^4u_6^4u_1^4u_2^4)$.
\end{theorem}

This evaluation, after a further specialisation $\kappa_i=(-1)^i\kappa$ gives an alternative to Theorem~\ref{Expression_P}.
Indeed, in the case $k=2\ell$, we look at a fixed $m\le\ell$ and collect all the coefficients of monomials $\kappa_{i_1}\dotsb\kappa_{i_{2m}}$ of length $2m$, which we further specialise $\kappa_i=(-1)^i$:
\begin{align*}
c_m
&=\sum_{1\le i_1<\dots<i_{2m}\le k}(-1)^{i_1+\dotsb+i_{2m}}
(1-x_{i_1}\dotsb x_{i_2-1})(1-x_{i_2}\dotsb x_{i_3-1})\dotsb
\\ &\qquad\qquad\qquad\qquad
\times(1-x_{i_{2m-1}}\dotsb x_{i_{2m}-1})(1-x_{i_{2m}}\dotsb x_{i_1-1}),
\end{align*}
where $x_j=u_j^4$ for $j=1,\dots,k$ as before.
Notice that $(-1)^{i_1+\dotsb+i_{2m}}$ is invariant under cyclic shifts.

\begin{lem} 
\label{Tri-diagonalDet}
For the determinant $D_k=D_k(u_1,u_2,\ldots, u_k)$ of the tridiagonal $(k-1)\times (k-1)$-matrix
\[
L_k= L_k(A_2^\circ,\dots,A_{k}^\circ; B_2,\ldots, B_{k-1}):=\begin{pmatrix}
A_2^\circ & B_2 & 0 & \hdots & 0 \\
B_2 & A_3 & B_3 & \hdots & 0 \\
0 & B_2 & A_4^\circ & \ddots & \vdots \\
\vdots & \vdots & \ddots & \ddots & B_{k-1} \\
0 & 0 & \hdots & B_{k-1} & A_{k}^\circ
\end{pmatrix}
\]
with $A_i^\circ= \frac{1-u_{i-1}^4u_i^4}{(1-u_{i-1}^4)(1-u_i^4)}$ \textup($i=2,\ldots, k$\textup)\textup, we have 
\[
D_k \cdot \prod_{j=1}^{k}(1-x_j)=1-x_1x_2\cdots x_k.
\]
\end{lem}

\begin{proof}
We prove the lemma by induction. It is obvious that $D_2\cdot (1-x_1)(1-x_2)= a_2\cdot (1-x_1)(1-x_2)= 1-x_1x_2$ and easy to see $D_3\cdot (1-x_1)(1-x_2)(1-x_3)=(A_2^\circ A_3^\circ-B_2^2)\cdot (1-x_1)(1-x_2)(1-x_3)=1-x_1x_2x_3$. 
Expanding the determinant $D_k$ along the last row leads to
\begin{align*}
D_k(u_1,u_2,\ldots, u_k)
& = A_k^\circ D_{k-1}(u_1,\ldots,u_{k-1})-B_{k-1}^2 D_{k-2}(u_1,\ldots,u_{k-2})
\\
& = \frac{1-x_{k-1}x_k}{(1-x_{k-1})(1-x_k)}D_{k-1}(u_1,\ldots, u_{k-1})-\frac{x_{k-1}}{(1-x_{k-1})^2} D_{k-2}(u_1,\ldots, u_{k-2}).
\end{align*}
Hence, we have 
\begin{align*}
D_k \cdot \prod_{j=1}^{k} (1-x_j)
&= \frac{1-x_{k-1}x_k}{1-x_{k-1}}D_{k-1}(u_1,\cdots,u_{k-1})\prod_{j=1}^{k} (1-x_j)
\\ &\quad 
- \frac{x_{k-1}}{1-x_{k-1}} (1-x_k) D_{k-2}(u_1,\ldots, u_{k-2}) \prod_{j=1}^{k-1} (1-x_j).
\end{align*}
It follows from the induction hypotheses 
\begin{align*}
D_k \cdot \prod_{j=1}^{k} (1-x_j)
&= \frac{1-x_{k-1}x_k}{1-x_{k-1}}(1-x_1\cdots x_{k-1}) -  \frac{x_{k-1}(1-x_k)}{1-x_{k-1}} (1-x_1\cdots x_{k-2})
\\
&= \frac1{1-x_{k-1}} 
\big(1-x_{k-1}x_k- x_1x_2\cdots x_{k-1}+ (x_1x_2\cdots x_k)\cdot x_{k-1}
\\ &\qquad
- x_{k-1}(1- x_k- x_1x_2\cdots x_{k-2} + (x_1\cdots x_{k-2}) \cdot x_k\big)
\\
&= 1- x_1x_2\cdots x_k,
\end{align*}
and the lemma follows. 
\end{proof}

\begin{proof}[Proof of Theorem~\textup{\ref{Expression_kappa}}]
Let us prove the theorem by showing explicitly the coefficients of $\kappa_{i_1}\cdots\kappa_{i_m}$ in $\det(R_k)\prod_{j=1}^k(1-x_j)$. 
First, remark that the coefficient of $\kappa_{i_1}\cdots\kappa_{i_m}$ is given by $\det(R_k)$ by putting $\kappa_j=0$ when $j\notin \{i_1, i_2,\ldots, i_m\}$. 
Note also that since the matrix is cyclic tridigaonal, we may assume that $i_1=1$. In fact, by cyclically moving the rows and columns in the determinant, we can easily see, for example, the following:
\[
\det \begin{pmatrix}
A_2 & B_2 & 0 & \hdots & B_1\\
B_2 & A_3 & B_3 & \hdots & 0 \\
0 & B_3 & A_4 & \ddots & \vdots \\
\vdots & \vdots & \ddots & \ddots & B_{k} \\
B_1 & 0 & \hdots & B_{k} & A_1
\end{pmatrix}
= \det \begin{pmatrix}
A_1 & B_1 & 0 & \hdots & B_k\\
B_1 & A_2 & B_2 & \hdots & 0 \\
0 & B_2 & A_3 & \ddots & \vdots \\
\vdots & \vdots & \ddots & \ddots & B_{k-1} \\
B_k & 0 & \hdots & B_{k-1} & A_k 
\end{pmatrix}.
\]
Then, for instance, the coefficient of $\kappa_1$ is obtained by  Lemma \ref{Tri-diagonalDet} as 
\[
\prod_{j=1}^k(1-x_j)\times \det L_{k}(A_2^\circ,\ldots,A_k^\circ; B_2,\ldots, B_{k-1})=1-x_1x_2x_3\cdots x_k.
\]
We consider two cases for the neighbor indices $i_{j-1}$ and $i_j$:
\begin{enumerate}
  \item[(i)] $i_j-i_{j-1}=1$
  \item[(ii)] $i_j-i_{j-1}>1$.
\end{enumerate}
For $j>1$ in the case (ii), we focus the following part of the matrix $M_k$:
\[
\begin{pmatrix}
A_{i_{j-1}}^\circ+\kappa_{i_{j-1}} & B_{i_{j-1}} & 0 & \hdots &0 & 0\\
B_{i_{j-1}}  & A_{i_{j-1}+1}^\circ +\kappa_{i_{j-1}+1}& B_{i_{j-1}+1} & \hdots &0 & 0 \\
0 & B_{i_{j-1}+1} & A_{i_{j-1}+2}^\circ +\kappa_{i_{j-1}+2} & \ddots &0 & \vdots \\
\vdots & \vdots & \ddots & \ddots & B_{i_j-2} &0 \\
0 & 0 & \hdots & B_{i_j-2} & A_{i_j-1}^\circ+\kappa_{i_j-1}& B_{i_j-1}\\
0 & 0 & \hdots & 0 & B_{i_j-1} & A_{i_j}^\circ+\kappa_{i_j}
\end{pmatrix}.
\]
We denote $L(i_{j-1}+1, i_j-1)$ the central part of the above matrix of size $(i_j-1)-(i_{j-1}+1)+1=i_j-i_{j-1} -1\,(\geq 1)$ by  putting $\kappa_{i_{j-1}+1}=\cdots= \kappa_{i_j-1}=0$. Then, it is elementary to see that 
the coefficient of  $\kappa_{i_1}\cdots\kappa_{i_m}$  in $\det(R_k)\prod_{j=1}^k(1-x_j)$ is given by 
\[
\prod_{j=1}^{m}\det L(i_{j-1}+1, i_j-1) \times \prod_{j=1}^k(1-x_j).
\] 
Here, we interpret the factor $\det L(i_{j-1}+1, i_j-1)=1$ for the cases (i), that is, when $i_j-i_{j-1}=1$. When $i_j-i_{j-1}>1$, by Lemma \ref{Tri-diagonalDet}, we  have 
$$
\det L(i_{j-1}+1, i_j-1) \times \prod_{n=i_{j-1}}^{i_j}(1-x_n)= 1- x_{i_{j-1}}x_{i_{j-1}+1}\cdots x_{i_j}.
$$
It follows that the coefficient of $\kappa_{i_1}\cdots\kappa_{i_m}$ in $\det(R_k)\prod_{j=1}^k(1-x_j)$ is given by 
\[
\prod_{j: i_j-i_{j-1}>1}(1- x_{i_{j-1}}x_{i_{j-1}+1}\cdots x_{i_j})\prod_{j: i_j-i_{j-1}=1}(1-x_{i_j}).
\]
This completes the proof of Theorem \ref{Expression_kappa}.
\end{proof}

\section{General $k\ge2\ell$}
\label{sec:general}

Define
\[
f_j(u_1,u_2,\dots,u_k)=f_j^{(k)}(u_1,u_2,\dots,u_k)=\prod_{i=1}^{j-1}u_i^2+\prod_{i=j}^ku_i^2,
\]
so that
$f_1(u_1,\dots,u_k)=1+(u_1\dotsb u_k)^2$, $f_2(u_1,\dots,u_k)=u_1^2+(u_2\dotsb u_k)^2$, \dots, $f_k(u_1,\dots,u_k)=(u_1\dotsb u_{k-1})^2+u_k^2$.

\begin{prop}
\label{th:Delta-inv}
The matrix $\Delta_k(\bu)$ has the following inverse:
{\small
\begin{align}
&
\Delta_k(\bu)^{-1}
=\frac1{1-(u_1\dotsb u_k)^2}
\nonumber\\ &\;\times
\begin{pmatrix}
f_1(u_1,u_2,\dots,u_k) & f_2(u_1,u_2,\dots,u_k) & f_3(u_1,u_2,\dots,u_k) & \hdots & f_k(u_1,u_2,\dots,u_k) \\
f_k(u_2,\dots,u_k,u_1) & f_1(u_2,\dots,u_k,u_1) & f_2(u_2,\dots,u_k,u_1) & \hdots & f_{k-1}(u_2,\dots,u_k,u_1) \\
f_{k-1}(u_3,\dots,u_1,u_2) & f_k(u_3,\dots,u_1,u_2) & f_1(u_3,\dots,u_1,u_2) & \hdots & f_{k-2}(u_3,\dots,u_1,u_2) \\
\vdots & \vdots & \vdots & \ddots & \vdots \\
f_2(u_k,u_1,\dots,u_{k-1}) & f_3(u_k,u_1,\dots,u_{k-1}) & f_4(u_k,u_1,\dots,u_{k-1}) & \hdots & f_1(u_k,u_1,\dots,u_{k-1})
\end{pmatrix},
\label{Delta-inv}
\end{align}}%
which can be compactly present as
\begin{equation*}
\Delta_k(\bu)^{-1}
=\bigg(\frac{f_{j-i+1}(u_i,u_{i+1},\dots,u_{i+k-1})}{1-(u_1\dotsb u_k)^2}\bigg)_{1\le i,j\le k}
\end{equation*}
after the cyclic index convention is implemented.
\end{prop}

\begin{proof}
Take $v_j=u_ju_{j+1}\dotsb u_k$ for $j=1,\dots,k$, so that $f_j(u_1,\dots,u_k)=v_1^2/v_j^2+v_j^2$, and denote by $M$ the $k\times k$ symmetric matrix on the right-hand side in \eqref{Delta-inv} without the pre-factor $1/(1-v_1^2)$.
Then the $(1,1)$ entry of the matrix $M\Delta_k(\bu)$ is equal to
\begin{align*}
&
f_1(u_1,\dots,u_k)\cdot\frac{1-u_k^4u_1^4}{(1-u_k^4)(1-u_1^4)}
-f_2(u_1,\dots,u_k)\cdot\frac{u_1^2}{1-u_1^4}
-f_k(u_1,\dots,u_k)\cdot\frac{u_k^2}{1-u_k^4}
\\ &\quad
=\frac{(1+v_1^2)(1-u_k^4u_1^4)}{(1-u_k^4)(1-u_1^4)}
-\frac{(u_1^2+v_1^2/u_1^2)u_1^2}{1-u_1^4}
-\frac{(v_1^2/u_k^2+u_k^2)u_k^2}{1-u_k^4}
%\\ &\quad
=1-v_1^2,
\end{align*}
while for $j\ge2$ the $(1,j)$ entry of $M\Delta_k(\bu)$ is
\begin{align*}
&
f_j(u_1,\dots,u_k)\cdot\frac{1-u_{j-1}^4u_j^4}{(1-u_{j-1}^4)(1-u_j^4)}
-f_{j-1}(u_1,\dots,u_k)\cdot\frac{u_{j-1}^2}{1-u_{j-1}^4}
-f_{j+1}(u_1,\dots,u_k)\cdot\frac{u_j^2}{1-u_j^4}
\\ &\quad
=\frac{(v_1^2/v_j^2+v_j^2)(1-u_{j-1}^4u_j^4)}{(1-u_{j-1}^4)(1-u_j^4)}
-\frac{(v_1^2/(u_{j-1}^2v_j^2)+u_{j-1}^2v_j^2)u_{j-1}^2}{1-u_{j-1}^4}
-\frac{(v_1^2u_j^2/v_j^2+v_j^2/u_j^2)u_j^2}{1-u_j^4}
%\\ &\quad
=0.
\end{align*}
The cyclic symmetry of both $M$ and $\Delta_k(\bu)$ implies that $M\Delta_k(\bu)=(1-v_1^2)\operatorname{Id}_k$, so that the result recorded in \eqref{Delta-inv} follows.
\end{proof}

We can use this explicit form from Proposition \ref{th:Delta-inv} for computing determinants
\begin{align*}
\det(\Delta_k(\bu)+\kappa\Xi)
&=\det\Delta_k(\bu)\cdot\det(\operatorname{Id}_k+\kappa\Xi\Delta_k(\bu)^{-1})
\\
&=\frac{(1-u_1^2\dotsb u_k^2)^2}{(1-u_1^4)\dotsb(1-u_k^4)}\cdot\det(\operatorname{Id}_k+\kappa\Xi\Delta_k(\bu)^{-1})
\end{align*}
when $\Xi$ is a diagonal matrix with $\pm\sqrt{-1}$ and $0$ along the diagonal.
In the case $\Xi=\Xi_k(i_1,i_2,\dots,i_{2\ell})$, the matrix product $\Xi\Delta_k(\bu)^{-1}$ happens to be a $k\times k$ matrix whose only nonzero rows are those numbered by $i_1<i_2<\dots<i_{2\ell}$; this implies that the determinant of $\operatorname{Id}_k+\kappa\Xi\Delta_k(\bu)^{-1}$ is equal to the determinant of the $(2\ell)\times(2\ell)$ matrix $\operatorname{Id}_{2\ell}+M$ where $M=(M_{r,s})_{1\le r,s\le2\ell}$ is the matrix with entries
\[
M_{r,s}=(-1)^r\sqrt{-1}\kappa\cdot\frac{f_{i_s-i_r+1}(u_{i_r},u_{i_r+1},\dots,u_{i_r+k-1})}{1-(u_1\dotsb u_k)^2}
\quad\text{for}\; 1\le r,s\le2\ell.
\]
Now observe that
\[
f_{i_s-i_r+1}^{(k)}(u_{i_r},u_{i_r+1},\dots,u_{i_r+k-1})
=f_{s-r+1}^{(2\ell)}(v_r,v_{r+1},\dots,v_{r+2\ell-1}),
\]
where $v_r=\prod_{j=i_r}^{i_{r+1}-1}u_j$, so that
\[
\det(\operatorname{Id}_{2\ell}+M)
=\det\big(\operatorname{Id}_{2\ell}+\kappa\Xi_{2\ell}(i_1,i_2,\dots,i_{2\ell})\Delta_{2\ell}(\bv)^{-1}\big).
\]
Since the latter determinant is computed previously, as the case $k=2\ell$, we obtain the same expressions \eqref{P2a} for the coefficients $c_m(\bu)$ in \eqref{poly-W} but with $\bv=(v_1,\dots,v_{2\ell})$ playing the role of $2\ell$~variables~$\bu$.
This leads to the following general formula for the coefficients.

\begin{theorem}
\label{th:gen}
For $k\ge2\ell$, we have
\begin{align*}
&
c_m(k;i_1,i_2,\dots,i_{2\ell};\bu)
=\frac12\sum_{r=1}^{m}\binom{\ell-r}{\ell-m}\sum_{j=r}^{2\ell-r}(-1)^{j-1}
\\ &\quad\times
\sum_{\substack{\bs:\len(\bs)=j\\\cohet(\bs)=r}}
\big(1-(v_{s_1}u_{s_2}\dotsb v_{s_j})^4\big)
\big(1-(u_1u_2\dotsb u_k)^4/(v_{s_1}u_{s_2}\dotsb v_{s_j})^4\big),
\end{align*}
where $v_s=\prod_{n=i_s}^{i_{s+1}-1}u_n$ for $s=1,2,\dots,2\ell$ and the range for $n$ in the latter product is understood cyclically.  
\end{theorem}

%\section{The leading part of a spectral zeta function}
%\label{sec:zeta}

\section{Concluding remarks}
\label{sec:conclusion}

Let us stress that the existence of integral expressions like \eqref{eq:Rkl} using polynomials $\mathcal{W}_k$ is somewhat unique for the special values of the NCHO spectral zeta function;
nothing similar is known for the quantum Rabi model mentioned in the Introduction.
Though both the NCHO and QRM are generalizations of the harmonic oscillator, the structure of their Hamiltonians is different, as those depend on two linearly independent matrices for the NCHO and on three matrices for the QRM.%
\footnote{To get integrals like  \eqref{eq:Rkl}, we need the techniques explored in \cite{IW05b} that only work for two linearly independent matrices.}
At the same time the NCHO represents a covering model of the QRM, and the recent work \cite{N26} brings to life the existence of series decompositions of type~\eqref{eq:specialvalues of zeta_Q(s)} for $\zeta_{\text{QRM}}(k)$ but with an \emph{infinite} sum on the right-hand side,
whose terms are derivatives of integrals that contain exponential rather than algebraic functions of $\bu$-like variables in the integrand.
Particular terms in the sum are already of great arithmetic interest:
Focusing on the first one (the first-order derivative) and applying a Taylor expansion to the exponential component, one finds that each resulting term can be written as a linear combination of generalizations of Beukers' integral from \cite{Beu1979}\,---\,the very integrals used in the study of the irrationality of $\zeta(2)$. Although this expression differs from the one involving $R_{k,\ell}(X)$, there are further similarities shared with our treatment of the polynomials $\mathcal{W}_k$ for the NCHO case. For instance, \cite[Corollary~4.4]{N26} features integrals reminiscent of the transfer-matrix trace examined in Section~\ref{sec:k=2l}.
A significant advancement can be to develop a direct link between the expressions of special values for the NCHO and (one- and two-photon) QRM cases using the similarities observed in~\cite{N26}.

\smallskip
The expressions like those in \eqref{eq:specialvalues of zeta_Q(s)}, \eqref{eq:Rkl} for the special values $\zeta_Q(s)$ may originate from our knowledge of the partition function of the NCHO.
Once an explicit formula for the latter is obtained, though it is still non-trivial (see, \cite[Conjecture 5.3]{KW25}), its Mellin's transform provides us with analytic continuation of $\zeta_Q(s)$ via contour integration\,---\,following familiar steps used for the Riemann zeta function $\zeta(s)$.
This leads to NCHO analogues of Bernoulli numbers and polynomials,
which can then be used to express special values at negative odd integers
(the values at even negative integers vanish).
Furthermore, these analogues are likely to be described in terms of hypergeometric functions. While one probably cannot expect a clean functional equation for $\zeta_Q(s)$ under the transformation $s\leftrightarrow 1-s$, it may still be possible to derive some results at positive integer points.

\smallskip
Though this is not transparent from the expressions for $R_{k,\ell}(X)$ in \eqref{eq:Rkl}, which in turn define $\zeta_{Q}(k)$ as given in~\eqref{eq:specialvalues of zeta_Q(s)}, there are further connections between the coefficients in the $\zeta(m,\frac12)$-expansions of $R_{2\ell,\ell}(X)$ and of $R_{2\ell+1,\ell}(X)$, more generally of $R_{k,\ell}(X)$ for $k>2\ell$. We make these reasonably explicit in \cite{WZ27} but also clarify their role in the modularity questions for $\zeta_{Q}(k)$ and in the parity-decomposition phenomenon.
We suspect\,---\,but do not have an explicit argument at the moment\,---\,that a reason behind such interconnections is the reduction mechanism exploited in the proof of Theorem~\ref{th:gen}.
Thus, a better understanding of the interrelation of the determinants $\mathcal{W}_k(X;u_1,\dots,u_k)$ when $k\ge2\ell$ varies and $\ell$ is fixed may shed further light on such connections.

\medskip
\noindent
\textbf{Acknowledgements.}
The authors thank the Max-Planck Institute for Mathematics (Bonn, Germany) for its hospitality and support during their visit there in June 2026. The second author further thanks the NTT Institute for Fundamental Mathematics (Japan) for its support during his visit in October 2026.

\medskip
\noindent
\textbf{AI involvement.}
Our original intention was to deduce Theorem~\ref{Expression_P} from Theorem~\ref{Expression_kappa} by the use of inclusion--exclusion principle.
The present simpler (and shorter) argument for proving Theorem~\ref{Expression_P}, namely the key identity of Proposition~\ref{key-prop} and the monomer--dimer matching expansion of a transfer-matrix trace, were offered to us by Claude Opus 5.5, even though it was not asked to do so but rather to verify numerically our calculations of the special zeta values.
At the end of the day we are thankful to this AI tool for its volunteering assistance.
The entire exposition including our proof of Theorem~\ref{Expression_P} was written without any AI help.

\bigskip

\begin{flushleft}
Masato Wakayama \par
NTT Institute of Fundamental Mathematics, \par 
NTT Inc.,\par 
3-9-11 Midori-cho Musashino-shi, Tokyo, 180-8585 JAPAN \par
%\texttt{masato.wakayama@ntt.com}

\medskip

Faculty of Social Informatics,\par
ZEN University, \par
4-12-15 Ginza Chuo-ku, Tokyo, 104-0061 JAPAN \par
\texttt{masato\_wakayama@zen.ac.jp}\par
\texttt{wakayama@imi.kyushu-u.ac.jp}

%\& \par
%Institute of Mathematics for Industry,\par
%Kyushu University \par
%744 Motooka, Nishi-ku, Fukuoka 819-0395 JAPAN \par
%\texttt{wakayama@imi.kyushu-u.ac.jp}

\bigskip

Wadim Zudilin \par
Institute for Mathematics, Astrophysics and Particle Physics, \par
 Radboud University Nijmegen, \par
PO Box 9010, 6500\,GL Nijmegen, NETHERLANDS \par
\texttt{w.zudilin@math.ru.nl} 

\bigskip

\end{flushleft}

\end{document}